\documentclass[11pt]{article}

\newcommand{\qed}{{\hfill\rule{4pt}{7pt}}}

\usepackage{amsmath, epsfig, cite}
\usepackage{amssymb,amsthm}
\usepackage{amsfonts}
\usepackage{latexsym}

\newtheorem{thm}{Theorem}[section]

\newtheorem{lem}[thm]{Lemma}

\numberwithin{equation}{section}

\makeatletter \@addtoreset{equation}{section} \makeatother

\begin{document}
\begin{center}
{\Large \bf Labeled Partitions and the $q$-Derangement Numbers}
\end{center}
\begin{center}
William Y. C. Chen$^{1}$ and  Deheng Xu$^{2}$

\vspace{3mm}
Center for Combinatorics, LPMC\\
Nankai University, Tianjin 300071\\
P. R. China

\small{Email: $^{1}$chen@nankai.edu.cn, $^{2}$xudeheng@eyou.com}
\end{center}

\noindent {\bf Abstract.} By a re-examination of  MacMahon's
original proof of his celebrated theorem on the distribution of the
major indices over permutations, we give a reformulation of his
argument in terms of the structure of labeled partitions. In this
framework, we are able to establish a decomposition theorem for
labeled partitions that leads to a simple bijective proof of
 Wachs' formula on the $q$-derangement numbers.


\noindent {\bf AMS Classification Numbers:} 05A30; 05A19, 05A15

\section{Introduction}

\newcommand{\maj}{{\rm maj}}

We will follow the terminology and notation on permutations and
partitions and $q$-series   in Andrews \cite{Andrews} and Stanley
\cite{Stan97}. The set of permutations on $\{1, 2, \ldots, n\}$ is
denoted by $S_n$. For any permutation $\pi=\pi_1\pi_2\cdots\pi_n\in
S_n$, an index $i$ with $1\le i \le n-1$ is called a \emph{descent}
of $\pi$ if $\pi_i>\pi_{i+1}$. The major index $\maj(\pi)$ of $\pi$,
introduced by MacMahon \cite{MacMahon},  is defined as the sum of
all descents of $\pi$. The following formula is well-known:
\begin{equation}
\label{e-MacMahon} \sum_{\pi\in S_n}q^{{\rm maj} (\pi)}=[n]!=1\cdot
(1+q) \cdot (1+q+q^2)\cdots (1+q+\cdots+q^{n-1}).
\end{equation}

The underlying idea of MacMahon' proof goes as follows. It is easier
to consider sequences and partitions  than solely permutations for
the purpose of studying the major index. MacMahon established
\eqref{e-MacMahon} by proving an equivalent formula
\begin{equation}\label{p-maj}
 \frac{1}{(q)_n}\sum_{\pi\in
S_n}q^{\maj(\pi)}=\frac{1}{(1-q)^n},
\end{equation}
where $(q)_n=(1-q)\cdots (1-q^{n-1})$, and $(q)_n^{-1}$ is the
generating function for partitions with at most $n$ parts. We will
give a reformulation of MacMahon's proof in Section 2 by introducing
the notion of standard labeled partitions.

The main objective of this paper is to employ MacMahon's method to
deal with the major index of derangements.  An integer $i$ with
$1\leq i\leq n$ is said to be a {\it fixed point} of $\pi\in S_n$ if
$\pi_i=i$, and {\it derangement point} otherwise.
\emph{Derangements} are permutations with no fixed points. Let $D_n$
be the set of all derangements in $S_n$. The $q$-derangement numbers
are defined by $d_0(q)=1$ and for $n\ge 1$
$$d_n(q)=\sum\limits_{\pi\in D_n}q^{\maj(\pi)}.$$

The following elegant formula was first derived by Gessel in his
manuscript and was published in \cite{GesselR} as a consequence of
the quasi-symmetric generating function encoding the descents and
the cycle structure of permutations. A combinatorial proof has been
obtained by Wachs \cite{wachs}:
\begin{equation}
\label{e-q-derangement}
d_n(q)=[n]!\sum\limits_{k=0}^n\frac{(-1)^k}{[k]!}q^{k\choose 2}.
\end{equation}

Let us review the combinatorial settings of Wachs for the above
formula.  Suppose the derangement points of $\pi$ are
$p_1,p_2,\cdots,p_{k}$. The reduction of $\pi$ to its
\emph{derangement part}, denoted by $dp(\pi)$, is defined as a
permutation on $\{1,2,\cdots,k\}$  induced by the relative order of
$\pi_{p_1},\pi_{p_2},\cdots,\pi_{p_k}$. For example, the derangement
points of $\pi=(1,5,3,7,6,2,9,8,4)$ are $2,4,5,6,7,9$, and
$\pi_2\pi_4\pi_5\pi_6\pi_7\pi_9=(5,7,6,2,9,4)$. Then
$dp(\pi)=(3,5,4,1,6,2)$. Clearly $dp(\pi)\in D_k$ if $\pi$ has $k$
derangement points. On the other hand, we can insert a fixed point
$j$ with $1\leq j\leq k+1$ into $\pi\in S_k$ to obtain a permutation
$$\bar{\pi}=\pi_1'\pi_2'\cdots\pi_{j-1}'j\ \pi_{j}'\cdots\pi_k'\in
S_{k+1},$$ where $\pi_i'=\pi_i$ if $\pi_i< j$ and $\pi_i'=\pi_i+1$
if $\pi_i\geq j$. Such an insertion operation produces one extra
fixed point.

 Wachs
\cite{wachs}  has established the following relation.

\begin{thm}\label{main}
Let $\sigma\in D_k$ and $k\leq n$. Then we have
\begin{equation}\label{target}
\sum\limits_{\substack{dp(\pi)=\sigma\\\pi\in
S_n}}q^{\maj(\pi)}=q^{\maj(\sigma)}{n\brack k},
\end{equation}
where ${n\brack k}=\frac{[n]!}{[k]![n-k]!}$ is the $q$-binomial
coefficient.
\end{thm}
By summing over all derangements $\sigma\in D_k$ and then summing
over all $k$ for the above relation \eqref{target}, and applying
\eqref{e-MacMahon} gives
$$[n]!=\sum\limits_{k=0}^n{n\brack k}d_k(q).$$
Thus \eqref{e-q-derangement} follows from the $q$-binomial
inversion \cite[Corollary 3.38]{Aigner},

In order to justify the relation \eqref{target},  Wachs found a
bijection on $S_n$ by rearranging a permutation $\pi $ according to
\emph{excedant} (where $\pi_i>i$), fixed point, and \emph{subcedant}
(where $\pi_i<i$). She showed that this bijection preserves the
major index by
considering $9$ cases. 
Then a result of Garsia-Gessel \cite[Theorem 3.1]{GG} on shuffles of
permutations is applied to establish
 Theorem \ref{main}.

Inspired by MacMahon's proof of \eqref{e-MacMahon}, we find it much
easier to deal with an equivalent form of \eqref{target}:
\begin{equation}\label{equbij}
\frac{1}{(q)_n}\sum\limits_{\substack{dp(\pi)=\sigma\\\pi\in
S_n}}q^{\maj(\pi)}=\frac{1}{(q)_k(q)_{n-k}}q^{\maj(\sigma)}.
\end{equation}

We will use the terminology of labeled partitions and will introduce
the notion of standard labeled partitions. In such terms, MacMahon's
proof can be easily stated. Moreover, a combinatorial reasoning of
\eqref{equbij} becomes quite natural, which is analogous to the
decomposition of a permutation by separating the derangements from
the fixed points.

\section{Labeled Partitions}

Let $\lambda=(\lambda_1,\dots,\lambda_n)$ be a partition, where
$\lambda_1\geq \lambda_2 \geq \cdots \geq \lambda_n\ge 0$. We say
that $\lambda$ is a partition with at most $n$ parts. We write
$|\lambda|=\lambda_1+\cdots +\lambda_n$. {\em A labeled partition}
is defined as a pair $(\lambda, \pi)$ of a partition $\lambda$ and a
permutation $\pi=\pi_1\pi_2\cdots \pi_n$. A labeled partition is
also represented in the following two row form as in Andrews
\cite[p. 43]{Andrews}:
\[ \left( \begin{array}{cccc}
      \lambda_1 & \lambda_2 & \cdots & \lambda_n \\
      \pi_1 & \pi_2 & \cdots & \pi_n
       \end{array}\right).
       \]
A labeled partition $(\lambda, \pi)$ is said to be \emph{standard}
if $\pi_{i}
> \pi_{i+1}$ implies $\lambda_i > \lambda_{i+1}$. For example, the
labeled partition in \eqref{exam-labeled-partition} is standard.

A labeled partition $(\lambda, \pi)$ is standard if $\lambda_i =
\lambda_{i+1}$ implies $\lambda_{i} < \lambda _{i+1}$.

The following Lemma \ref{lemma}  is straightforward to verify, which
is MacMahon's approach to study the major index with the aid of
partitions, see MacMahon \cite{MacMahon}, Andrews \cite[Theorem
3.7]{Andrews}, Knuth \cite[p.~18]{Knuth1} or \cite{Knuth}. This
method was further extended by Stanley \cite{Stan72}. For other
applications, see \cite{GG}.

\begin{lem}\label{lemma}
Given $\pi\in S_n$, there is a bijection $\psi_{\pi}\colon
\lambda\mapsto\mu$ from partitions $\lambda$ with at most $n$ parts
to standard labeled partitions $(\mu, \pi)$ such that
$|\lambda|+\maj(\pi)=|\mu|$.
\end{lem}

The bijection $\psi_\pi$ (or simply $\psi$ when $\pi$ is understood
from the context), is given as follows:
$$\mu=\psi_\pi(\lambda)=(\lambda_1+\phi_1,\lambda_2+\phi_2,\cdots,\lambda_n+\phi_n),$$
where $\phi_i$ is the number of descents in
$\pi_i\pi_{i+1}\cdots\pi_n$. One may also view  $\psi$
as the operation of adding 1 to
$\lambda_1,\dots, \lambda_i$ whenever $i$ is a descent of
$\pi$.

We now give a restatement of MacMahon's proof of (\ref{p-maj}) in
the above terminology of labeled partitions.

\noindent
{\it Proof of $(\ref{p-maj})$.} Given a
sequence $a_1 a_2\cdots a_n$ of nonnegative integers,
we associate it with a weight $q^{a_1+a_2+\cdots+a_n}$.
Let us construct a two row array
\[ \left( \begin{array}{cccc}
      a_1 & a_2 & \cdots & a_n \\
       1 &  2 & \cdots & n
       \end{array}\right).
       \]
By permuting the columns of the above array, one can get a unique
standard labeled partition $(\mu, \pi)$ with $|\mu|=a_1+a_2+\cdots
+a_n$.
 Applying Lemma \ref{lemma},
we obtain a partition $\lambda$ with $|\lambda|+\maj(\pi) =\mu$.
Clearly, the above steps are reversible. This completes the proof.
\qed

\noindent{\bf An Example.} Let $a_1a_2\dots a_9$ be the sequence
with a two line array
\[\left(
\begin{array}{ccccccccc}
       3 & 6 &   8 & 3 & 1 & 3 & 6 & 4 & 8 \\
       1 & 2 &   3 & 4 & 5 & 6 & 7& 8 & 9
       \end{array}\right).\]
Permuting the columns we get the a standard labeled partition:
\begin{equation}
{\mu\choose \pi}= \left( \begin{array}{ccccccccc}
       8 & 8 &   6 & 6 & 4 & 3 & 3 & 3 & 1 \\
       3 & \underline{9} &   2 & 7 & \underline{8} & 1 & 4 & \underline{6} & 5
       \end{array}\right), \label{exam-labeled-partition}
\end{equation}
where we have underlined the descents of $\pi$.

Applying $\psi^{-1}$ gives
\[{\lambda\choose \pi}=\left( \begin{array}{ccccccccc}
       5 & 5 &   4 & 4 & 2 & 2 & 2 & 2 & 1 \\
       3 & \underline{9} &   2 & 7 & \underline{8} & 1 & 4 & \underline{6} & 5
       \end{array}\right).\]

\medskip
We remark that the idea of standard labeled partitions appeared in
\cite[p.~292]{GG}, though it was not used to prove (\ref{p-maj}).

We now come to the main result of this note, which is a
decomposition theorem on standard labeled partitions in terms of the
fixed points. Let ${\mu \choose \pi}$ be a standard labeled
partition with $\pi\in S_n$.  Assume that $\pi$ has $n-k$ fixed
points. Let $i_1< i_2< \cdots < i_{n-k}$ be the fixed points, let
$j_1 < j_2 < \cdots< j_k$ be the derangement points of $\pi$, and
let $dp(\pi)=\sigma \in D_k$. We now define the following
decomposition of a standard labeled partition:
\begin{equation}
   \varphi \colon  {\mu \choose \pi} \mapsto (\beta,\gamma),
\end{equation}
where $\beta=\mu_{j_1}\mu_{j_2}\cdots \mu_{j_k}$ and
$\gamma=\mu_{i_1} \mu_{i_2} \cdots \mu_{i_{n-k}}$ are the partitions
corresponding to derangement points and fixed points, respectively.
Evidently $|\mu|=|\beta|+|\gamma|$.

The following is the main theorem of this paper.

\begin{thm}\label{p-main} For given $\sigma \in D_k$, the
decomposition $\varphi$ of ${\mu\choose \pi}$ with $dp(\pi)=\sigma$
is a bijection from standard labeled partitions to pairs of
partitions such that $(\beta, \sigma)$ is a standard labeled
partition.
\end{thm}

We note that the above theorem and Lemma \ref{lemma} lead to a
combinatorial interpretation of the relation (\ref{equbij}). Since
$(\beta, \sigma)$ is a standard labeled partition, we may find a
partition $\alpha$ such that $\psi{\alpha \choose \sigma} = {\beta
\choose \sigma}$.  Consequently, the bijection $\varphi\circ \psi$
maps a labeled partition ${\lambda \choose \pi}$ to a pair $(\alpha,
\gamma)$ of partitions, where $\alpha$ has at most $k$ parts and
$\gamma$ has at most $n-k$ parts. Moreover, the following relation
holds:
\begin{equation}
\lambda+\maj(\pi) =  |\alpha|+|\gamma| +\maj(\sigma),
\end{equation}
which implies (\ref{equbij}).

\noindent {\it Proof of Theorem \ref{p-main}.} We first show that
$(\beta,\sigma)$ is standard. It suffices to show that if
$\pi_i>\pi_j$ with $\pi_{i+1},\dots,\pi_{j-1}$ being fixed points,
then $\mu_i>\mu_j$. If $j=i+1$, since ${\mu \choose \pi}$ is
standard, we have $\mu_i > \mu_j$. We now consider the case $i <
j-1$, and we claim that either $\pi_i>\pi_{i+1}=i+1$ or
$\pi_{j-1}=j-1>\pi_j$ holds; Otherwise, it follows that $\pi< i+1
\leq j-1< \pi_j$, a contradiction. Therefore, we have either $\mu_i>
\mu_{i+1}$ or $\mu_{j-1} > \mu_j$. It is deduced that $\mu_i>
\mu_j$.

We now proceed to construct the map $\varphi'$ which is guided by
the procedure of inserting the fixed points of $\pi$ to the
derangement $\sigma$ on $\{1, 2, \ldots, k\}$. We will show that
$\varphi'$ and $\varphi$ are inverse to each other, which implies
that $\varphi$ is a bijection.

 Let
$(\mu^0,\pi^0)=(\beta,\sigma)$. We assume that $(\mu^{i},\pi^{i})$
is obtained from $(\mu^{i-1},\pi^{i-1})$ by \emph{inserting}
$\gamma_i$. We find the first position $r$ so that the insertion of
$\gamma_{i}$ at the proper position produces a partition. This
partition is denoted $\mu^{i}$. Clearly,
$\mu^i_{r-1}>\mu^i_r=\gamma_{i}$. Assume that
$\mu^i_r=\cdots=\mu^i_{t}>\mu^i_{t+1}$ for some $t\ge r$. If $r=t$
then we set $s=r$. Otherwise we find the position $s$ such that
$\pi^{i-1}_{s-1}<s\le \pi_{s}^{i-1}$, (here we have taken
${\pi^{i-1}_{r-1}}$ as $-\infty$ and ${\pi^{i-1}_{t}}$ as $\infty$).
Now insert $s$ as a fixed point into $\pi^{i-1}$ to generate
$\pi^{i}$. Note that the position $s$ is judiciously chosen so that
the subsequence $\pi^i_r,\pi^i_{r+1},\cdots,\pi^i_{t}$, which is the
same as ${\pi^{i-1}_r}',\cdots ,{\pi^{i-1}_{s-1}}', s,
{\pi_{s}^{i-1}}',\cdots ,{\pi_{t-1}^{i-1}}'$, is increasing, and
hence $\pi^i$ is a standard labeled partition.

Since $\mu^{n-k}$ is the partition obtained from $\beta$ by
inserting $\gamma_1,\dots,\gamma_{n-k}$, we must have
$\mu^{n-k}=\mu$. From the above procedure, one sees that $\pi^{n-k}$
is constructed from $\pi_0=\sigma$ by inserting fixed points,
therefore we have $dp(\pi^{n-k})=\sigma$. It follows that for a
given $\sigma\in D_k$, we have
$\varphi\varphi'(\beta,\gamma)=(\beta,\gamma)$.

Now it is only necessary to show that $\pi^{n-k}=\pi$. For
simplicity, we write $\pi^{n-k}$ as $\bar{\pi}$. We prove by
contradiction. By removing same fixed points, we may assume that the
first fixed point $f$ of $\pi$ is different from the first fixed
point $\bar{f}$ of $\bar{\pi}$. Furthermore, we may assume that
$f<\bar{f}$. Clearly, $\mu_{f}=\mu_{\bar{f}}$. Since $(\mu,\pi)$ and
$(\mu,\bar{\pi})$ are standard labeled partitions, we have
$$\pi_{f}<\pi_{f+1}<\cdots<\pi_{\bar{f}},\textrm{ and } \bar{\pi}_{f}<\bar{\pi}_{f+1}<\cdots<\bar{\pi}_{\bar{f}}.$$
So we get $\bar{\pi}_{f}=\sigma_f \ge \pi_{f+1}-1\geq \pi_{f}=f$. By
assumption, $f$ is not a fixed point of $\bar{\pi}$. It follows that
$\bar{\pi}_f>f$. Hence $\bar{\pi}_{\bar{f}}>\bar{f}$, a
contradiction. \qed

\noindent {\bf An Example.}

Let $\displaystyle{\lambda \choose \pi}={5\ 4\ 4\ 4\ 4\ 3\
2\choose\underline{5}\ \underline{2}\ 1\ 4\ \underline{7}\ 3\ 6}.$
Applying $\psi$, we get $\displaystyle{\mu \choose \pi}={8\ 6\ 5\ 5\
5\ 3\ 2\choose5\ 2\ 1\ 4\ 7\ 3\ 6}$.

The fixed points of $\pi$ are $2,4$. Hence $\sigma=dp(\pi)=(3\ 1\ 5\
2\ 4)$.  Applying $\varphi$ on $(\mu,\pi)$ gives
$(\beta,\gamma)=((8\ 5\ 5\ 3\ 2),(6\ 5)).$ Finally, applying
$\psi^{-1}$ to $(\beta,\sigma)$, we obtain $ {\alpha\choose
\sigma}={6\ 4\ 4\ 3\ 2\choose3\ 1\ 5\ 2\ 4}.$ Based on $\sigma=(3\
1\ 5\ 2\ 4)$, we conclude that $\varphi{\lambda \choose \pi}=((6\ 4\
4\ 3\ 2),\ (6\ 5))$.

\medskip
Conversely, given $\sigma=(3\ 1\ 5\ 2\ 4)$ and $(\beta,\gamma)=((8\
5\ 5\ 3\ 2),\ (6\ 5))$, we have
$${\beta\choose \sigma}={8\ 5\ 5\ 3\ 2\choose3\ 1\ 5\ 2\
4}\stackrel{\gamma_1=6}{\longmapsto}{8\ \underline{6}\ 5\ 5\ 3\
2\choose4\ \underline{2}\ 1\ 6\ 3\
5}\stackrel{\gamma_2=5}{\longmapsto}{8\ 6\ 5\ \underline{5}\ 5\ 3\
2\choose5\ 2\ 1\ \underline{4}\ 7\ 3\ 6}={\mu\choose \pi}.$$

\vskip 6mm

\noindent {\bf Acknowledgments.} We would like to thank Guoce Xin
for valuable suggestions. This work was supported by the 973 Project
on Mathematical Mechanization, the Ministry of Education, the
Ministry of Science and Technology and the National Science
Foundation of China.

\end{document}